\documentclass[11pt,reqno]{amsart}
\usepackage{amsmath,amsthm,amssymb,verbatim,enumerate,ifthen,times}
\usepackage[mathscr]{eucal}
\usepackage[utf8]{inputenc}
\usepackage[T1]{fontenc}
\newcommand{\N}{\mathbb{N}}

\newcommand{\R}{\mathbb{R}}

\newcommand{\B}{\mathscr{B}}

\newcommand{\Cdot}{\!\cdot\!}

\DeclareMathOperator{\sgn}{sign}

\newtheorem{theorem}{Theorem}
\newtheorem*{theorem*}{Theorem}
\def\Thm#1#2{\ifthenelse{\equal{#1}{*}}{\begin{theorem*}#2\end{theorem*}}
  {\begin{theorem}\label{T#1}#2\end{theorem}}}
\newtheorem{Atheorem}{Theorem}

\def\thm#1{Theorem~\ref{T#1}}

\newtheorem{proposition}[theorem]{Proposition}
\newtheorem*{proposition*}{Proposition}
\def\Prp#1#2{\ifthenelse{\equal{#1}{*}}{\begin{proposition*}#2\end{proposition*}}
             {\begin{proposition}\label{P#1}#2\end{proposition}}}

\newtheorem{corollary}[theorem]{Corollary}
\newtheorem*{corollary*}{Corollary}
\def\Cor#1#2{\ifthenelse{\equal{#1}{*}}{\begin{corollary*}#2\end{corollary*}}
             {\begin{corollary}\label{C#1}#2\end{corollary}}}
\def\cor#1{Corollary~\ref{C#1}}

\newtheorem{lemma}[theorem]{Lemma}
\newtheorem*{lemma*}{Lemma}
\def\Lem#1#2{\ifthenelse{\equal{#1}{*}}{\begin{lemma*}#2\end{lemma*}}
             {\begin{lemma}\label{L#1}#2\end{lemma}}}
\def\lem#1{Lemma~\ref{L#1}}

\newtheorem{example}[theorem]{Example}
\newtheorem*{example*}{Example}
\def\Exa#1#2{\ifthenelse{\equal{#1}{*}}{\begin{example*}\rm #2\end{example*}}
             {\begin{example}\label{Ex#1}\rm #2\end{example}}}

\newtheorem{problem}[theorem]{Problem}

\theoremstyle{definition}
\newtheorem{definition}[theorem]{Definition}

\newtheorem{remark}[theorem]{Remark}
\newtheorem*{remark*}{Remark}
\def\Rem#1#2{\ifthenelse{\equal{#1}{*}}{\begin{remark*}\rm #2\end{remark*}}
             {\begin{remark}\label{R#1}\rm #2\end{remark}}}

\newcommand{\eq}[1]{\eqref{E#1}}
\newcommand{\Eq}[2]{\ifthenelse{\equal{#1}{*}}
  {\begin{equation*}\begin{aligned}[]#2\end{aligned}\end{equation*}}
  {\begin{equation}\begin{aligned}[]\label{E#1}#2\end{aligned}\end{equation}}}

\long\def\comment#1{}

\begin{document}
\large

\date{\today}

\title[Equality of Bajraktarevi\'c means to quasi-arithmetic means]
{On the equality of Bajraktarevi\'c means to quasi-arithmetic means}

\author[Zs. P\'ales]{Zsolt P\'ales}
\address{Institute of Mathematics, University of Debrecen,
H-4002 Debrecen, Pf.\ 400, Hungary}
\email{pales@science.unideb.hu}

\author[A. Zakaria]{Amr Zakaria}
\address{Department of Mathematics, Faculty of Education, Ain Shams University, Cairo 11341, Egypt; Institute of Mathematics, University of Debrecen, H-4002 Debrecen, Pf.\ 400, Hungary}
\email{amr.zakaria@edu.asu.edu.eg}



\thanks{The research of the first author was supported by the Hungarian Scientific Research Fund (OTKA) Grant K-111651 and by the EFOP-3.6.1-16-2016-00022 project. This project is co-financed by the European Union and the European Social Fund.}

\begin{abstract}
This paper offers a solution of the functional equation 
\Eq{*}{
 \big(tf(x)+(1-t)f(y)\big)\varphi(tx+(1-t)y)=tf(x)\varphi(x)+(1-t)f(y)\varphi(y)
 \qquad(x,y\in I),
}
where $t\in\,]0,1[\,$, $\varphi:I\to\R$ is strictly monotone, and $f:I\to\R$ is an arbitrary unknown function. As an immediate application, we shed new light on the equality problem of Bajraktarevi\'c means with quasi-arithmetic means. 
\end{abstract}
\maketitle

\section{Introduction}

Throughout this paper, the symbols  $\R$, $\R_+$, and $\N$ will stand for the sets of real, positive real, and natural numbers, respectively, and $I$ will always denote a nonempty open real interval.

For $n\in\N$, define the set of $n$-dimensional weight vectors $\Lambda_n$ by
\Eq{*}{
  \Lambda_n:=\{(\lambda_1,\dots,\lambda_n)\in \R^n\mid\lambda_1,\dots,\lambda_n\geq0,\,\lambda_1+\dots+\lambda_n>0\}.
}
A function $M:I^n\times\Lambda_n\to I$ is called an \emph{$n$-variable weighted mean} if, for all $x=(x_1,\dots,x_n)\in I^n$ and $\lambda=(\lambda_1,\dots,\lambda_n)\in \Lambda_n$,
\Eq{*}{
  \min\big\{x_i\mid \lambda_i>0\big\}
  \leq M(x,\lambda)\leq\max\big\{x_i\mid \lambda_i>0\big\}.
}
The most classical class of weighted means is the class of power means, or more generally, quasi-arithmetic means. Their definition is recalled from the book \cite{HarLitPol34}.

Given a continuous strictly monotone function $\varphi:I\to\R$, the \emph{weighted  quasi-arithmetic mean} $A_\varphi:\bigcup_{n=1}^\infty I^n\times\Lambda_n\to I$ is defined by
\Eq{*}{
  A_\varphi(x,\lambda)
  :=\varphi^{-1}\left(\frac{\lambda_1\varphi(x_1)+\dots+ \lambda_n\varphi(x_n)}{\lambda_1+\dots+\lambda_n}\right)
}
for $n\in\N$, $x=(x_1,\dots,x_n)\in I^n$, and $\lambda=(\lambda_1,\dots,\lambda_n)\in\Lambda_n$. The restriction of $A_\varphi$ to the set $I^n\times\Lambda_n$ is called the \emph{$n$-variable weighted  quasi-arithmetic mean}. In the case when $\lambda_1=\dots=\lambda_n=1 $, we speak about an \emph{$n$-variable (discrete) quasi-arithmetic mean} and write $A_\varphi(x)$ instead of $A_\varphi(x,\lambda)$. The function $\varphi$ is called the \emph{generating function} of the quasi-arithmetic mean $A_\varphi$.

By taking $\varphi(x):=x$ for $x\in\R$, the resulting mean $A_\varphi$ is the \emph{weighted arithmetic mean}. Given $p\in\R$, $p\neq0$, the function $\varphi(x):=x^p$ $(x\in\R_+)$ generates the \emph{$p$th weighted power mean}. To obtain the \emph{weighted geometric mean}, one should take the weighted quasi-arithmetic mean generated by $\varphi(x):=\log(x)$ $(x\in\R_+)$.

For the equality of quasi-arithmetic means, we have the following equivalence of six conditions.

\Thm{1}{(\cite{HarLitPol34}, \cite{Kuc85})
Let $\varphi,\psi:I\to\R$ be continuous strictly monotone functions.
Then the following properties are pairwise equivalent:
\begin{enumerate}[(i)]
 \item $A_\varphi(x,\lambda)=A_\psi(x,\lambda)$ holds for all $n\geq2$, $x=(x_1,\dots,x_n)\in I^n$ and $\lambda=(\lambda_1,\dots,\lambda_n)\in\Lambda_n$.
 \item $A_\varphi(x)=A_\psi(x)$ for all $n\geq2$ and $x=(x_1,\dots,x_n)\in I^n$.
 \item $A_\varphi(x,\lambda)=A_\psi(x,\lambda)$ holds for all $x=(x_1,x_2)\in I^2$ and $\lambda=(\lambda_1,\lambda_2)\in\Lambda_2$.
 \item $A_\varphi(x)=A_\psi(x)$ holds for all $x=(x_1,x_2)\in I^2$.
 \item There exists $t\in\,]0,1[\,$, such that $A_\varphi(x,\lambda)=A_\psi(x,\lambda)$ holds for all $x=(x_1,x_2)\in I^2$ with $\lambda=(t,1-t)\in\Lambda_2$.
 \item There exist $a,b\in\R$ such that $\psi=a\varphi+b$.
\end{enumerate}}

Generalizing the notion of quasi-arithmetic means, Mahmud Bajraktarevi\'c in 1958 introduced a new class of means in the following way: Let $\varphi:I\to\R$ be a continuous strictly monotone function, let $f:I\to\R_+$ be a positive function and define $A_{\varphi,f}:\bigcup_{n=1}^\infty I^n\times\Lambda_n\to I$ by
\Eq{*}{
  A_{\varphi,f}(x,\lambda)
  :=\varphi^{-1}\left(\frac{\lambda_1f(x_1)\varphi(x_1)+\dots+\lambda_nf(x_n)\varphi(x_n)}{\lambda_1f(x_1)+\dots+\lambda_nf(x_n)}\right)
}
for $n\in\N$, $x=(x_1,\dots,x_n)\in I^n$, and $\lambda=(\lambda_1,\dots,\lambda_n)\in\Lambda_n$. Due to the identity
\Eq{*}{
  A_{\varphi,f}((x_1,\dots,x_n),(\lambda_1,\dots,\lambda_n))
  =A_{\varphi}((x_1,\dots,x_n),(\lambda_1 f(x_1),\dots,\lambda_n f(x_n)),
}
one can immediately see that the restriction of the function $A_{\varphi,f}(x,\lambda)$ to the set $I^n\times\Lambda_n$ is an $n$-variable weighted mean.

Denoting $g:=\varphi\cdot f$, we can rewrite $A_{\varphi,f}(x,\lambda)$ in the following more symmetric form:
\Eq{*}{
  B_{g,f}(x,\lambda)
  :=\Big(\frac{g}{f}\Big)^{-1}\left(\frac{\lambda_1g(x_1)+\dots+\lambda_ng(x_n)}{\lambda_1f(x_1)+\dots+\lambda_nf(x_n)}\right).
}
In fact, if $g$ is also nowhere zero, then one can see that $B_{g,f}\equiv B_{f,g}$. It is also clear that the expression for $B_{g,f}$ is well defined if $f$ is positive and $g/f$ is strictly monotone and continuous.

In order to describe necessary and sufficient conditions for the equality of Bajraktarevi\'c means, we introduce the following terminology. We say that two pairs of functions $(f,g):I\to\R^2$ and $(h,k):I\to\R^2$ are equivalent (and we write $(f,g)\sim(h,k)$) if there exist constants $a,b,c,d$ with $ad\neq cb$ such that
\Eq{abcd}{
  h=af+bg \qquad\mbox{and}\qquad k=cf+dg.
}
One can easily check that $\sim$ is an equivalence relation, indeed. 

For two given functions $f,g:I\to\R$, we define the two-variable function $\Delta_{f,g}:I^2\to\R$ as follows
\Eq{*}{
  \Delta_{f,g}(x,y):=\left|\begin{array}{cc}
  f(x) & f(y) \\ 
  g(x) & g(y)
  \end{array}\right| \qquad(x,y\in I).
}

For the equality of Bajraktarevi\'c means, we have the following equivalence of four conditions.

\Thm{BM}{(\cite{AczDar63c}, \cite{DarPal82})
Let $f,g,h,k:I\to\R$ such that $f$ and $h$ are positive functions and $g/f$, $k/h$ are continuous and strictly monotone.
Then the following properties are pairwise equivalent:
\begin{enumerate}[(I)]
 \item $B_{g,f}(x,\lambda)=B_{k,h}(x,\lambda)$ holds for all $n\geq2$, $x=(x_1,\dots,x_n)\in I^n$ and $\lambda=(\lambda_1,\dots,\lambda_n)\in\Lambda_n$.
 \item $B_{g,f}(x)=B_{k,h}(x)$ for all $n\geq2$ and $x=(x_1,\dots,x_n)\in I^n$.
 \item $B_{g,f}(x,\lambda)=B_{k,h}(x,\lambda)$ holds for all $x=(x_1,x_2)\in I^2$ and $\lambda=(\lambda_1,\lambda_2)\in\Lambda_2$.
 \item[(VI)] $(f,g)\sim(h,k)$.
\end{enumerate}}

The proof of the above theorem is partly based on the following lemma that we will also need in the sequel.

\Lem{EP}{Let $(f,g):I\to\R^2$ and $(h,k):I\to\R^2$ be equivalent pairs. Then, for some nonzero constant $\gamma$,  
\Eq{d2}{
  \Delta_{h,k}=\gamma\Delta_{f,g}.
}}

\begin{proof} By the assumption, there exist constants $a,b,c,d$ with $ad\neq cb$ such that \eq{abcd} holds. Then, using the product theorem of determinants, for all $x,y\in I$,
\Eq{*}{
  \Delta_{h,k}(x,y)
  =\left|\begin{array}{cc}
  h(x) & h(y) \\ 
  k(x) & k(y)
  \end{array}\right|
  =\left|\begin{array}{cc}
   a & b \\ 
   c & d
  \end{array}\right| \Cdot 
  \left|\begin{array}{cc}
  f(x) & f(y) \\ 
  g(x) & g(y)
  \end{array}\right|
  =(ad-bc)\Delta_{f,g}(x,y).
}
Therefore, \eq{d2} holds with $\gamma:=ad-bc\neq0$.
\end{proof}

When comparing the characterizations of the equality for quasi-arithmetic and Bajraktarevi\'c means, one can observe that two conditions are missing from the list of \thm{BM} (which would correspond to assertions (iv) and (v) in \thm{1}):  
\begin{enumerate} \it
 \item[(IV)] $B_{g,f}(x)=B_{k,h}(x)$ holds for all $x=(x_1,x_2)\in I^2$.
 \item[(V)] there exists $t\in\,]0,1[\,$, such that $B_{g,f}(x,\lambda)=B_{k,h}(x,\lambda)$ holds for all $x=(x_1,x_2)\in I^2$ with $\lambda=(t,1-t)\in\Lambda_2$.
\end{enumerate}
It is obvious that each of the equivalent assertions (I), or (II), or (III), or (VI) implies (IV). It is also evident that (IV) implies (V) (with $t:=\frac12$). As it has been pointed out in our paper \cite{PalZak19b}, assertion (V) with $t\in\,]0,\frac12[\,\cup\,]\frac12,1[\,$ implies (VI) (and hence also (I) and (II) and (III)) under three times differentiability of the generating functions $f,g,h$, and $k$. On the other hand, as it was shown by Losonczi \cite{Los99}, assertion (IV) is not equivalent to any of the assertions (I), (II), (III), and (VI). More precisely, under six times differentiability, Losonczi completely described the solutions of the equality problem of two-variable Bajraktarevi\'c means and established 32 cases of the equality beyond the standard equivalence of the generating pairs.

Similar problems have been considered in the literature by several authors. Bajraktarevi\'c \cite{Baj58}, \cite{Baj69} solved the equality problem of two Bajraktarevi\'c means with at least three variables under three times differentiability. He also found sufficient conditions for the equality of the two-variable means. Acz\'el and Dar\'oczy \cite{AczDar63c} described the necessary and sufficient conditions of the equality for all number of variables but without imposing any additional regularity properties. Dar\'oczy and Losonczi \cite{DarLos70} solved the comparison problem assuming first-order differentiability. Losonczi \cite{Los99} solved the equality problem of two-variable Bajraktarevi\'c assuming a certain algebraic conditions and six times differentiability of the unknown functions. Later, he \cite{Los06b} investigated the equality problem of more general means under the same regularity assumptions, but he removed the algebraic conditions required in his earlier papers. In a recent paper by Losonczi and P\'ales \cite{LosPal11a}, the equality of two-variable Bajraktarevi\'c means generated via two different measures has been investigated. 
Until now, the weakening of the regularity assumptions has not been succeeded in the general case, only in the particular case when the equality problem of (symmetric) two-variable Bajraktarevi\'c mean with a quasi-arithmetic mean was considered. Matkowski \cite{Mat12a} 2012 solved this question supposing first-order differentiability. He did not notice however, that the same goal was accomplished 8 years ago in 2004 by Daróczy, Maksa and Páles \cite{DarMakPal04} where no additional differentiability condition was assumed. 

The goal of this paper to solve the above mentioned equality problem in a particular case but without additional unnatural regularity assumptions. More precisely, we will solve the equality problem of Bajraktarevi\'c means to quasi-arithmetic means in two settings: in the class of two-variable symmetric means and in the class of two-variable nonsymmetrically weighted or more than three-variable weighted means. After an obvious substitution, these equality problems can be reduced to the functional equation
\Eq{FE}{
 \big(tf(x)+(1-t)f(y)\big)\varphi(tx+(1-t)y)=tf(x)\varphi(x)+(1-t)f(y)\varphi(y)
 \qquad(x,y\in I),
}
where $f,\varphi:I\to\R$ and $t\in\,]0,1[\,$ is fixed. This equation was considered and solved in the case $t=\frac12$ in \cite{DarMakPal04} under strict monotonicity and continuity of $\varphi$ and in \cite{KisPal18} under continuity 
of $\varphi$, respectively. In \thm{0} and \thm{0+} below, we completely solve \eq{FE} assuming only the strict monotonicity of $\varphi$ and also including the case $t\neq\frac12$. Applying these solutions, the main results are stated in \thm{M3} and \thm{M4}, which provide various equivalent conditions for a Bajraktarevi\'c mean to be quasi-arithmetic.

\section{Solution of the fundamental functional equation \eq{FE}}

\Thm{0}{Let $\varphi:I\to\R$ be a strictly monotone function, $f:I\to\R$ be an arbitrary function, and $t\in\,]0,1[\,$. Assume that the functional equation \eq{FE} holds. Then either $f$ is identically zero, or $f$ is nowhere zero, $f$ and $\varphi$ are infinitely many times differentiable and there exists a nonzero constant $\gamma\in\R$ such that
\Eq{Sol}{
  f^2\varphi'=\gamma.
}}

\begin{proof} If $f$ is identically zero, then \eq{FE} holds, therefore no information can be obtained for $\varphi$.

Assume now that there exists a point $y_0$ such that $f$ does not vanish at $y_0$. Then, for $x\in I$, $x\neq y_0$, the convex combination $tx+(1-t)y_0$ is strictly between the values $x$ and $y_0$. Therefore, by the strict monotonicity of $\varphi$, we have that $(\varphi(tx+(1-t)y_0)-\varphi(x))(\varphi(y_0)-\varphi(tx+(1-t)y_0))>0$. Then, it follows from the functional equation \eq{FE}, that
\Eq{fx}{
  f(x)=\frac{1-t}t f(y_0)\frac{\varphi(y_0)-\varphi(tx+(1-t)y_0)}
  {\varphi(tx+(1-t)y_0)-\varphi(x)}.
}
This implies that $f(x)$ is nonzero for all $x\in I$, furthermore, $f(x)$ has the same sign as $f(y_0)$, i.e., the sign of $f$ is constant. 

In what follows, we prove that, at every point of $I$, the function $f$ has left and right limits and it is continuous at every point where $\varphi$ is continuous. Denote by $D_\varphi$ the set of discontinuity points of $\varphi$. Then the monotonicity of $\varphi$ implies that $D_\varphi$ is countable.

Let $x_0\in I$ be fixed. Then $tx_0+(1-t)I$ is a subinterval of $I$, hence $I\setminus D_\varphi$ intersects $tx_0+(1-t)I$. Therefore, there exists an element $y_0\in I$ such that $tx_0+(1-t)y_0\in I\setminus D_\varphi$. Thus, $\varphi$ is continuous at $tx_0+(1-t)y_0$. Now, upon taking the left or right limits as $x$ tends to $x_0$ of the right hand side of equality \eq{fx}, we can see that these limits exist because $\varphi(tx+(1-t)y_0)$ tends to $\varphi(tx_0+(1-t)y_0)$ and $\varphi(x)$ has a left and right limit (by the monotonicity of $\varphi$).  Therefore, \eq{fx} yields that $f$ has left and right limits at $x_0$. In addition, if $\varphi$ is continuous at $x_0$, then its left and right limits are the same, hence $f$ has to be continuous at $x_0$. 

From what we have proved it follows that $f$ is continuous everywhere except at countably many points, hence $f$ is continuous almost everywhere. On the other hand, $f$ is bounded on every compact subinterval of $I$. Indeed, if $f$ were unbounded on a compact subinterval $[a,b]\subseteq I$, then there would exist a subsequence $(x_n)$ in $[a,b]$ converging to some element $x_0\in[a,b]$, such that $|f(x_n)|\to +\infty$. We can extract a subsequence $(x_{n_k})$ which is either converging from the left or from the right to $x_0$. Then the limit of $f(x_{n_k})$ is the left or right limit of $f$ at $x_0$, which is finite, contradicting $|f(x_{n_k})|\to +\infty$. Having the local boundedness of $f$, it follows that $f$ is Riemann integrable on every compact subinterval of $I$.

Let $0<\alpha<\frac12|I|$ and $I_\alpha:=(I-\alpha)\cap(I+\alpha)$. Then $I_\alpha$ is a nonempty interval and $I_\alpha+[-\alpha,\alpha]\subseteq I$. Let $u\in I_\alpha$, $v\in[-\alpha,\alpha]$ and substituting $x:=u-(1-t)v$ and $y:=u+tv$ into \eq{FE}, we obtain that
\Eq{*}{
 \big(tf\big(u-(1-t)v\big)&+(1-t)f\big(u+tv\big)\big)\varphi(u)\\
 &=tf\big(u-(1-t)v\big)\varphi\big(u-(1-t)v\big)+(1-t)f\big(u+tv\big)\varphi\big(u+tv\big).
}
holds for all $u\in I_\alpha$ and for all $v\in[-\alpha,\alpha]$. Integrating both sides of the previous equation on $v\in[0,\alpha]$
it follows that
\Eq{*}{
 \varphi(u)\int_{-\alpha}^{\alpha}&\big(tf\big(u-(1-t)v\big)+(1-t)f\big(u+tv\big)\big)dv\\
 &=t\int_{-\alpha}^{\alpha}f\big(u-(1-t)v\big)\varphi\big(u-(1-t)v\big)dv+(1-t)\int_{-\alpha}^{\alpha}f\big(u+tv\big)\varphi\big(u+tv\big)dv.
}
After simple change of the variable transformations, for all $u\in I_\alpha$, we get
\Eq{phiu}{
 \varphi(u)\bigg(-\frac{t}{1-t}\int_{u-(1-t)\alpha}^{u+(1-t)\alpha}f
 +\frac{1-t}{t}\int_{u-t\alpha}^{u+t\alpha}f\bigg)
 &=-\frac{t}{1-t}\int_{u-(1-t)\alpha}^{u+(1-t)\alpha}f\cdot\varphi
 +\frac{1-t}{t}\int_{u-t\alpha}^{u+t\alpha}f\cdot\varphi.
}
Having that $f$ is either positive everywhere or negative everywhere, it follows that $\varphi(u)$ is the ratio of two expressions that are
continuous with respect to $u$. Therefore, $\varphi$ and hence $f$ are continuous everywhere in $I_\alpha$. This, together with \eq{phiu}, implies that $\varphi(u)$ is the ratio of two expressions that are
continuously differentiable with respect to $u$. Hence $\varphi$ is continuously differentiable on $I_\alpha$. Since $0<\alpha<\frac12|I|$ is arbitrary, it follows that $\varphi$ is continuously differentiable and $f$ is continuous on $\bigcup_{\alpha>0} I_\alpha=I$. Going back to formula \eq{fx}, the continuous differentiability of $\varphi$ implies that $f$ is also continuously differentiable.

Now, we show that $\varphi$ and $f$ are twice continuously differentiable. Differentiating \eq{FE} with respect to $x$, we have
\Eq{FE'}{
f'(x)\varphi(tx+(1-t)y)+(tf(x)+(1-t)f(y))\varphi'(tx+(1-t)y)=(f\varphi)'(x)\qquad(x,y\in I).
} 
By substituting $x:=u-(1-t)v$ and $y:=u+tv$ into the previous equation and integrating both sides on $v\in[0,\alpha]$, we get
\Eq{*}{
\varphi'(u)\int_{-\alpha}^{\alpha}&\big(tf(u-(1-t)v)+(1-t)f(u+tv)\big)dv\\
&=-\varphi(u)\int_{-\alpha}^{\alpha}f'(u-(1-t)v)dv+\int_{-\alpha}^{\alpha}(f\varphi)'(u-(1-t)v)dv\qquad(u\in I_\alpha).
} 
After similar change of the variable transformations as \eq{phiu}, for all $u\in I_\alpha$, we obtain
\Eq{*}{
 \varphi'(u)\bigg(-\frac{t}{1-t}\int_{u-(1-t)\alpha}^{u+(1-t)\alpha}f
 +\frac{1-t}{t}\int_{u-t\alpha}^{u+t\alpha}f\bigg)
 &=\frac{1}{1-t}\varphi(u)\int_{u-(1-t)\alpha}^{u+(1-t)\alpha}f'
 -\frac{1}{1-t}\int_{u-(1-t)\alpha}^{u+(1-t)\alpha}(f\varphi)'.
}
From here it follows that $\varphi'$ is the ratio of two continuously differentiable functions on $I_\alpha$. Thus $\varphi$ is twice continuously differentiable on $I_\alpha$ and hence on $I$. This result, combined with \eq{fx}, implies that $f$ is two times continuously differentiable on $I$.

To prove that $\varphi$ and $f$ are infinitely many times differentiable, differentiate \eq{FE'} with respect to $y$, to get
\Eq{2oc}{
  (f'(x)+f'(y))\varphi'(tx+(1-t)y)+(tf(x)+(1-t)f(y))\varphi''(tx+(1-t)y)=0.
}
Substituting $y:=x$, we arrive at 
\Eq{2oc+}{
2f'\varphi'+f\varphi''=0,
}
or equivalently,
\Eq{*}{
(f^2\varphi')'=0.
}
Hence there exists a real constant $\gamma$ such that $f^2\varphi'=\gamma$. If $\gamma$ were zero, then this equation would imply that $\varphi'$ is identically zero, which contradicts the strict monotonicity of $\varphi$. As a consequence,  \eq{Sol} holds. Finally, applying \eq{Sol} and \eq{fx} repeatedly, we get that $\varphi$ and $f$ are infinitely many times differentiable. 
\end{proof}

In order to describe the solution of functional equation \eq{FE}, we introduce the following notation. 

For a real parameter $p\in\R$, introduce the sine and cosine type functions $S_p,C_p:\R\to\R$ by
\Eq{*}{
  S_p(x):=\begin{cases}
           \sin(\sqrt{-p}x) & \mbox{ if } p<0, \\
           x & \mbox{ if } p=0, \\
           \sinh(\sqrt{p}x) & \mbox{ if } p>0, 
         \end{cases}\qquad\mbox{and}\qquad
  C_p(x):=\begin{cases}
           \cos(\sqrt{-p}x) & \mbox{ if } p<0, \\
           1 & \mbox{ if } p=0, \\
           \cosh(\sqrt{p}x) & \mbox{ if } p>0. 
        \end{cases}
}
It is easily seen that the functions $S_p$ and $C_p$ form the fundamental system of solutions for the second-order homogeneous linear differential equation $h''=ph$.

\Thm{0+}{Let $\varphi:I\to\R$ be a strictly monotone function, $f:I\to\R$ be a non-identically-zero function, and $t\in\,]0,1[\,$. Then the following assertions are equivalent:
\begin{enumerate}[(i)]
 \item $(\varphi,f)$ solves \eq{FE};
 \item $f$ is nowhere zero, $f$ and $\varphi$ are twice differentiable such that \eq{2oc+} holds and there exists $p\in\R$ with $(t-\frac12)p=0$ such that $f''=pf$;
 \item $f$ is nowhere zero and there exists $p\in\R$ with $(t-\frac12)p=0$ such that 
\Eq{SC}{
(f,f\Cdot\varphi)\sim(S_p,C_p).
}
\end{enumerate}}

\begin{proof}
Assume that $(\varphi,f)$ solves \eq{FE}. Then, as we have proved in \thm{0}, our conditions imply that $f$ is nowhere zero, $f$ and $\varphi$ are infinitely many times differentiable, and there exists a nonzero $\gamma\in\R$ such that \eq{Sol} holds. As in the proof of \thm{0}, differentiating \eq{FE} with respect to $x$ and then with respect to $y$, we get equations \eq{FE'} and \eq{2oc}, respectively. Substituting $y:=x$ into the last equality, \eq{2oc+} follows immediately.

Differentiating \eq{2oc} with respect to $x$, we obtain
\Eq{3oc}{
f''(x)\varphi'(tx+(1-t)y)&+(2tf'(x)+tf'(y))\varphi''(tx+(1-t)y)\\
&+t(tf(x)+(1-t)f(y))\varphi'''(tx+(1-t)y)=0.
}
Inserting $y:=x$, it follows that
\Eq{*}{
f''\varphi'+t\big(3f'\varphi''+f\varphi'''\big)=0.
}
On the other hand, differentiating \eq{2oc+} with respect $x$, we obtain
\Eq{3oc++}{
  2f''\varphi'+3f'\varphi''+f\varphi'''=0.
}
Combining the above equalities, we conclude that 
\Eq{1/2}{
(1-2t)f''\varphi'=0.
}
Due to \eq{Sol}, $\varphi'$ is nowhere zero. Consequently, either $t=\frac12$ or $f''=0$ on $I$.

In the first case when $t\neq\frac12$, then $f''=0$, and hence, assertion (ii) holds with $p=0$.

In the case $t=\frac12$, equation \eq{1/2} does not provide any information on $f$ and $\varphi$. Therefore, we substitute $t=\frac12$ into \eq{3oc}, to get
\Eq{*}{
f''(x)\varphi'\Big(\dfrac{x+y}{2}\Big)+(f'(x)+\frac12 f'(y))\varphi''\Big(\dfrac{x+y}{2}\Big)
+\frac14(f(x)+f(y))\varphi'''\Big(\dfrac{x+y}{2}\Big)=0.
}
Differentiating this equality with respect to $y$, we obtain
\Eq{*}{
\frac12(f''(x)+f''(y))\varphi''\Big(\dfrac{x+y}{2}\Big)+\frac12(f'(x)+f'(y))\varphi'''\Big(\dfrac{x+y}{2}\Big)+\frac18 (f(x)+f(y))\varphi''''\Big(\dfrac{x+y}{2}\Big)=0.
}
Substituting $y:=x$ and multiplying by $4$, we arrive at 
\Eq{4oc}{
4f''\varphi''+4f'\varphi'''+ f\varphi''''=0.
}
However, differentiating \eq{3oc++}, we obtain
\Eq{*}{
2f'''\varphi'+5f''\varphi''+4f'\varphi'''+f\varphi''''=0.
}
Subtracting \eq{4oc} from this equality side by side, we get 
\Eq{*}{
2f'''\varphi'+f''\varphi''=0.
}
Using \eq{Sol} and \eq{2oc+}, we can eliminate $\varphi'$ and $\varphi''$, and thus we get
\Eq{*}{
\dfrac{f'''f-f''f'}{f^2}=0.
}
Equivalently,
\Eq{*}{
\Big(\dfrac{f''}{f}\Big)'=0,
}
which implies that there exists a constant $p\in\R$ such that $f''=pf$.
This proves the last part of statement (ii).

Assume now that assertion (ii) holds, i.e., $f$ is nowhere zero, equation \eq{2oc+} and $f''=pf$ hold for some constant $p\in\R$ with $(t-\frac12)p=0$. Therefore, there exist constants $a,b\in\R$ such that
\Eq{fun2}{
  f=a S_p+b C_p
}
On the other hand, using equation \eq{2oc+}, it follows that
\Eq{*}{
  (f\Cdot\varphi)''
  =f''\Cdot\varphi+2f'\Cdot\varphi'+f\Cdot\varphi''
  =f''\Cdot\varphi=p f\Cdot\varphi,
}
which means that $g:=f\Cdot\varphi$ satisfies the differential equation $g''=pg$. Hence, there exist constants $c,d\in\R$ such that
\Eq{fun3}{
  f\Cdot\varphi= cS_p+dC_p.
}
From the two equalities \eq{fun2} and \eq{fun3}, it follows that $(f,f\Cdot\varphi)\sim(S_p,C_p)$, that is, assertion (iii) holds.

Finally, assume that (iii) is valid. Then $f$ is nowhere zero on $I$ and the equivalence \eq{SC} holds on $I$ for some $p\in\R$ with $(t-\frac12)p=0$. This, in view of \lem{EP}, implies that there exists a nonzero constant $\gamma$ such that
\Eq{*}{
  \Delta_{f,f\cdot\varphi}=\gamma\Delta_{S_p,C_p}.
}
On the other hand, the functional equation \eq{FE} holds if and only if 
\Eq{*}{
 t\Delta_{f,f\cdot\varphi}(x,tx+(1-t)y)
 +(1-t)\Delta_{f,f\cdot\varphi}(y,tx+(1-t)y)=0 \qquad(x,y\in I).
}
Therefore, to complete the proof, it is sufficient to prove that 
\Eq{d3}{
 t\Delta_{S_p,C_p}(x,tx+(1-t)y)
 +(1-t)\Delta_{S_p,C_p}(y,tx+(1-t)y)=0 \qquad(x,y\in I).
}
In the case $p=0$, we have that
\Eq{*}{
 t\Delta_{S_0,C_0}(x,tx+(1-t)y)+(1-t)\Delta_{S_0,C_0}(y,tx+(1-t)y)
 =\left|\begin{array}{cc}
  tx+(1-t)y & tx+(1-t)y \\ 
  1& 1
  \end{array}\right|=0.
}
In the case $t=\frac12$ and $p<0$, denote $q:=\sqrt{-p}$. Using well-known identities for trigonometric functions, we get
\Eq{*}{
  &\frac12\Delta_{S_p,C_p}\Big(x,\dfrac{x+y}{2}\Big)
 +\frac12\Delta_{S_p,C_p}\Big(y,\dfrac{x+y}{2}\Big)\\
 &=\left|\begin{array}{cc}
  \dfrac{\sin(qx)+\sin(qy)}2 & \sin\Big(q\dfrac{x+y}{2}\Big) \\ 
   \dfrac{\cos(qx)+\cos(qy)}2& 
   \cos\Big(q\dfrac{x+y}{2}\Big)
  \end{array}\right|
  =\left|\begin{array}{cc}
  \sin\Big(q\dfrac{x+y}{2}\Big)\cos\Big(q\dfrac{x-y}{2}\Big) 
  & \sin\Big(q\dfrac{x+y}{2}\Big) \\ 
   \cos\Big(q\dfrac{x+y}{2}\Big)\cos\Big(q\dfrac{x-y}{2}\Big)& 
   \cos\Big(q\dfrac{x+y}{2}\Big)
  \end{array}\right|=0.
}
Similar arguments apply to the case $p>0$ by using identities for hyperbolic functions, and therefore we leave it to the reader to verify \eq{d3}.
\end{proof}

Given an at most second-degree polynomial $P(u):=\alpha +\beta u+\gamma u^2$, where $\alpha,\beta,\gamma\in\R$, we call the value $D_P:=\beta^2-4\alpha\gamma$ the \emph{discriminant of $P$}. 

\Lem{DP}{If $P$ is an at most second-degree polynomial, then  $D_P=(P')^2-2P''P$.}

\begin{proof} Let $P$ be of the form $P(u):=\alpha +\beta u+\gamma u^2$, where $\alpha,\beta,\gamma\in\R$. Then
\Eq{*}{
  ((P')^2-2P''P)(u)=(\beta+2\gamma u)^2-4\gamma(\alpha +\beta u+\gamma u^2)
  =\beta^2-4\alpha\gamma=D_P,
}
which was to be proved.
\end{proof}

The following result is instrumental for our main results. 

\Lem{P}{Let $P$ be an at most second degree polynomial which is positive on $I$ let $D_P$ denote its discriminant and let $t\in\,]0,1[\,$ with $(t-\frac12)D_P=0$. Let $\psi$ be a primitive function of $1/P$ and $\ell:=1/\sqrt{P}$. Then the functions $\varphi:=\psi^{-1}$ and $f:=\ell\circ\varphi$ satisfy equation \eq{FE} on the interval $\psi(I)$.}

\begin{proof}
In order to prove that $(\varphi,f)$ solves \eq{FE}, we show that \thm{0+} part (ii) is valid. An easy computation shows that
\Eq{P}{
\varphi'=\frac{1}{\psi'\circ\psi^{-1}}=\frac{1}{\psi'\circ\varphi},
\qquad\varphi''=-\frac{\psi''\circ\psi^{-1}}{(\psi'\circ\psi^{-1})^3}=-\frac{\psi''\circ\varphi}{(\psi'\circ\varphi)^3},\qquad\mbox{and}\qquad 
\psi'=\frac{1}{P}=\ell^2. 
}
Therefore, it is obvious that 
\Eq{*}{
f^2\cdot\varphi'=(\ell^2\circ\varphi)\cdot\frac{1}{\psi'\circ\varphi}=1.
}
As a consequence, after differentiating both sides, we get that \eq{2oc+} holds.
Now, we only need to show that there exists $p\in\R$ such that $(t-\frac12)p=0$ and $f''=pf$. After simple calculations, using \eq{P} and \lem{DP}, we get 
\Eq{*}{
f''&=(\ell\circ\varphi)''=(\ell''\circ\varphi)\varphi'^2+(\ell'\circ\varphi)\varphi''=\bigg(\ell''\cdot\frac{1}{(\psi')^2}-\ell'\cdot\frac{\psi''}{(\psi')^3}\bigg)\circ\varphi \\
&=\bigg(\dfrac{\ell''\ell-2\ell'^2}{\ell^5}\bigg)\circ\varphi
=\bigg(\frac{P'^2-2P''P}{4\sqrt{P}}\bigg)\circ\varphi
=\bigg(\frac{D_P}{4\sqrt{P}}\bigg)\circ\varphi=\frac{D_P}{4}f.
}
Consequently, with $p:=D_P/4$ the equality $f''=pf$ holds on $\psi(I)$ and hence assertion (ii) of \thm{0+} is satisfied. 
\end{proof}

\section{Main results}

For simplicity, we introduce the following regularity classes for the generating functions of Bajraktarevi\'c means as follows: Let the class
$\B_0(I)$ contain all pairs $(f,g)$ such that
\begin{enumerate}
 \item[(i)] $f$ is everywhere positive on $I$.
 \item[(ii)] $g/f$ is strictly monotone and continuous on $I$.
\end{enumerate}
For $n\geq1$, let $\B_n(I)$ denote the class of all pairs $(f,g)$ such that
\begin{enumerate}
 \item[(+i)] $f$ is everywhere positive on $I$ and $f,g:I\to\R$ are $n$ times continuously differentiable functions.
 \item[(+ii)] $(g/f)'$ is nowhere zero on $I$.
\end{enumerate}
For $(f,g)\in\B_n(I)$ and for $i,j\in\{0,\dots,n\}$, we introduce the following notations:
\Eq{W}{
  W^{i,j}_{f,g}:=\left|\begin{array}{cc} f^{(i)} & f^{(j)} \\ g^{(i)} & g^{(j)} \end{array}\right|
  \qquad\mbox{and}\qquad
  \Phi_{f,g}:=\frac{W^{2,0}_{f,g}}{W^{1,0}_{f,g}},\qquad
  \Psi_{f,g}:=-\frac{W^{2,1}_{f,g}}{W^{1,0}_{f,g}}.
}
The following lemma was stated and verified in \cite{PalZak19a}.

\Lem{DE}{Let $(f,g)\in\B_2(I)$. Then $f,g$ form a fundamental system of solutions of the second-order homogeneous linear differential equation
\Eq{DE}{
  Y''=\Phi_{f,g}Y'+\Psi_{f,g}Y.
}}

As a consequence of \thm{0+}, we can immediately obtain a characterization of the equality between two-variable weighted Bajraktarevi\'c means and  
two-variable weighted quasi-arithmetic means.

\Cor{M3}{
Let $t\in\,]0,1[\,$, $(f,g)\in\B_0(I)$, and let $h:I\to\R$ be a continuous strictly monotone function. Then
\Eq{B}{
  B_{g,f}((x,y),(t,1-t))=A_h((x,y),(t,1-t)) \qquad (x,y\in I)
}
holds if and only if there exists $p\in\R$ with $(t-\frac12)p=0$ such that
\Eq{gf}{
  (f,g)\sim(S_p\circ h,C_p\circ h).
}
}

\begin{proof}
Applying $g/f$ to the both sides of \eq{B} and substituting $F:=f\circ h^{-1}$, $G:=g\circ h^{-1}$, and $\varphi:=\dfrac{G}{F}$, we get an equivalent formulation of \eq{B} as follows:
\Eq{B+}{
 \big(tF(u)+(1-t)F(v)\big)\varphi(tu+(1-t)v)=tF(u)\varphi(u)+(1-t)F(v)\varphi(v)
 \qquad(u,v\in h(I)).
}
Thus, the pair $(\varphi,F)$ satisfies \eq{FE} on the interval $h(I)$.
Therefore, by \thm{0+}, $p\in\R$ with $(t-\frac12)p=0$ such that $(F,G)=(F,F\Cdot\varphi)\sim(S_p,C_p)$ holds on $h(I)$. After substitution, this yields that \eq{gf} holds on the interval $I$.
\end{proof}

The last two theorems contain the main results of our paper. They offer various characterizations of the equality of a Bajraktarevi\'c mean to a quasi-arithmetic mean. In the first result we consider such an equality for the (symmetric) two-variable setting.

\Thm{M3}{
Let $(f,g)\in\B_0(I)$. Then the following statements are equivalent.
\begin{enumerate}[(i)]
\item There exists a continuous strictly monotone function $h:I\to\R$ such that
\Eq{fgh}{
  B_{g,f}(x,y)=A_h(x,y)\qquad(x,y\in I).
}
\item There exist real constants $\alpha,\beta,\gamma$ such that
\Eq{abc}{
  \alpha f^2+\beta fg+\gamma g^2=1.
}
\item Provided that $(f,g)\in\B_1(I)$, equation \eq{fgh}
holds with $h=\int W^{1,0}_{f,g}$.
\item Provided that $(f,g)\in\B_2(I)$, there exists a real constant $\delta$ such that
\Eq{W++}{ 
  W^{2,1}_{f,g}=\delta (W^{1,0}_{f,g})^3.
}
\item Provided that $(f,g)\in\B_2(I)$, $\Psi_{f,g}$ is differentiable and
\Eq{Psi}{
\Psi'_{f,g}=2\Phi_{f,g}\Psi_{f,g}.
}
\end{enumerate}
}

\begin{proof} We will prove first the equivalence of statements (i) and (ii).

Assume first that (i) holds, i.e., there exists a continuous strictly monotone function $h:I\to\R$ such that \eq{fgh} is valid. Then, applying \cor{M3} for $t=\frac12$, it follows that there exists $p\in\R$ such that the equivalence in \eq{gf} holds. Therefore, there exist $a,b,c,d\in\R$ with $ad\neq bc$ such that
\Eq{abcd+}{
  S_p\circ h=af+bg \qquad\mbox{and}\qquad C_p\circ h=cf+dg.
}
Using well-known trigonometric and hyperbolic identities, we have that
\Eq{*}{
  C_p^2-\sgn(p)\Cdot S_p^2=1
}
holds on $\R$, and hence $C_p^2\circ h-\sgn(p)\Cdot S_p^2\circ h=1$ holds on $I$. Combining this identity with the equalities in \eq{abcd+}, we get
\Eq{*}{
  (cf+dg)^2-\sgn(p)\Cdot (af+bg)^2=1
}
on $I$. Therefore, statement (ii) holds with 
\Eq{*}{
  \alpha:=c^2-\sgn(p)a^2,\qquad \beta:=2cd-2\sgn(p)\Cdot ab, \qquad\mbox{and}\qquad  \gamma:=d^2-\sgn(p)b^2.
}

Assume now that assertion (ii) is valid, i.e., \eq{abc} holds with some real  constants $\alpha,\beta,\gamma$. Denote $\varphi:=g/f$. Then, by $(f,g)\in\B_0(I)$, we have that $\varphi$ is strictly monotone and continuous. Replacing $g$ by $f\Cdot\varphi$ in \eq{abc}, we get
\Eq{PP0}{
  \alpha+\beta\varphi+\gamma\varphi^2=\frac1{f^2}.
}
Hence
\Eq{PP}{
  P(u):=\alpha+\beta u+\gamma u^2=\frac1{f^2\circ\varphi^{-1}(u)} \qquad(u\in\varphi(I)).
}
Thus, $P$ is an at most second-degree polynomial which is positive on the interval $J:=\varphi(I)$. Now, we are in the position to apply \lem{P} in the case $t=\frac12$. 
Let $\psi$ be a primitive function of $1/P$ and $\ell:=1/\sqrt{P}$. Then the functions $\varphi^*:=\psi^{-1}$ and $f^*:=\ell\circ\varphi^*$ satisfy equation \eq{FE} on $\psi(J)$. This immediately implies that the two-variable Bajraktarevi\'c mean $B_{f^*\Cdot\varphi^*,f^*}$ equals the two-variable arithmetic mean on $\psi(J)$, that is, for all $u,v\in\psi(J)$,
\Eq{*}{
  (\varphi^*)^{-1}\bigg(\frac{f^*(u)\varphi^*(u)+f^*(v)\varphi^*(v)}
  {f^*(u)+f^*(v)}\bigg)=\frac{u+v}{2}.
}
Now substituting $u:=(\varphi^*)^{-1}\circ\varphi(x)$ and $v:=(\varphi^*)^{-1}\circ\varphi(y)$ where $x,y\in I$, and observing that
\Eq{*}{
  f^*\circ (\varphi^*)^{-1}\circ\varphi
  =\ell\circ\varphi=\frac{1}{\sqrt{P\circ\varphi}}
  =\sqrt{f^2\circ\varphi^{-1}\circ\varphi}=f,
}
the above equality, for all $x,y\in I$, implies that
\Eq{*}{
  (\varphi^*)^{-1}\bigg(\frac{f(x)\varphi(x)+f(y)\varphi(y)}{f(x)+f(y)}\bigg)
  =\frac{(\varphi^*)^{-1}\circ\varphi(x)+(\varphi^*)^{-1}\circ\varphi(y)}{2}.
}
Applying the function $\varphi^{-1}\circ\varphi^*$ to this equation side by side, it follows that the two-variable Bajraktarevi\'c mean $B_{f\Cdot\varphi,f}$ equals the two-variable quasi-arithmetic mean $A_h$ on $I$, where $h:=(\varphi^*)^{-1}\circ\varphi$.

The implication (iii)$\Rightarrow$(i)$\sim$(ii) is obvious. Therefore, it remains to prove the implication (ii)$\Rightarrow$(iii). Assume that $(f,g)\in\B_1(I)$.
If (ii) holds for some $\alpha,\beta,\gamma\in\R$, then define the polynomial $P$ by \eq{PP} and let $\psi:=\int(1/P)$. As we have seen it before, then (i) holds with $h:=-(\varphi^*)^{-1}\circ\varphi=-\psi\circ\varphi$. Therefore,
\Eq{*}{
  h'=-(\psi'\circ\varphi)\Cdot\varphi'
  =-\bigg(\frac{1}{P}\circ\varphi\bigg)\Cdot\bigg(\frac{g}{f}\bigg)'
  =\frac{W_{f,g}^{1,0}}{(P\circ\varphi)\Cdot f^2}=W_{f,g}^{1,0}.
}
This completes the proof of assertion (iii)

To prove the implication (ii)$\Rightarrow$(iv), assume that $(f,g)\in\B_2(I)$.
If (ii) holds for some $\alpha,\beta,\gamma\in\R$, then equation \eq{PP0} is satisfied, where $P$ is the polynomial defined in \eq{PP} and hence $\frac1{f^2}=P\circ\varphi$. Differentiating this equality once and twice, it follows that
\Eq{*}{
  -2\frac{f'}{f^3}=(P'\circ\varphi)\Cdot\varphi' \qquad\mbox{and}\qquad
  \frac{6(f')^2-2ff''}{f^4}
  =(P''\circ\varphi)\Cdot(\varphi')^2+(P'\circ\varphi)\Cdot\varphi''.
}
Solving this system of equations with respect to $P'\circ\varphi$ and $P''\circ\varphi$, we obtain
\Eq{*}{
  P'\circ\varphi=-2\frac{f'}{f^3\varphi'} \qquad\mbox{and}\qquad
  P''\circ\varphi=\frac{6(f')^2\varphi'-2ff''\varphi'+2ff'\varphi''}{f^4(\varphi')^3}.
}
On the other hand, we have the following two equalities
\Eq{*}{
  \varphi'=\bigg(\frac{g}{f}\bigg)'=\frac{fg'-f'g}{f^2}
  =-\frac{W_{f,g}^{1,0}}{f^2}
}
and
\Eq{*}{
W^{2,1}_{f,g}=W^{2,1}_{f,f\cdot\varphi}
=\left|\begin{array}{cc} f'' & f' \\ (f\varphi)'' & (f\varphi)' \end{array}\right|
=-2(f')^2\varphi'+ff''\varphi'-ff'\varphi''.
}
Therefore, using \lem{DP}, we get
\Eq{*}{
  D_P&=(P'\circ\varphi)^2-2(P''\circ\varphi)(P\circ\varphi)
  =\frac{-8(f')^2\varphi'+4ff''\varphi'-4ff'\varphi''}{(f^2\varphi')^3}
  =\frac{4W_{f,g}^{2,1}}{\big(-W_{f,g}^{1,0}\big)^3},
}
which shows that (iv) holds with $\delta:=-D_P/4$. 

To prove the implication (iv)$\Rightarrow$(i), let $(f,g)\in\B_2(I)$. If (iv) holds for some real constant $\delta$, then 
\Eq{Psi++}{
\Psi_{f,g}=-\delta (W^{1,0}_{f,g})^2.
}
Let $Y\in\{f,g\}$. Then, as we have stated it in \lem{DE}, $Y$ is a solution of the second-order homogeneous linear differential equation \eq{DE}. In view of \eq{Psi++}, this differential equation is now of the form
\Eq{DE+}{
  Y''=\Phi_{f,g} Y'-\delta (W^{1,0}_{f,g})^2Y.
}
In order to solve this equation, let $\xi$ be an arbitrarily fixed point of the interval $I$, define $h:I\to\R$ by $h(x):=\int_\xi^x W^{1,0}_{f,g}$. Then $h$ is twice differentiable and strictly monotone with a nonvanishing first derivative, hence its inverse is also twice differentiable. Now define $Z:=Y\circ h^{-1}$. Then $Z:h(I)\to\R$ is a twice differentiable function and we have $Y=Z\circ h$. Differentiating $Y$ once and twice, we get
\Eq{*}{
Y'=(Z'\circ h)h'\qquad\mbox{and}\qquad Y''=(Z''\circ h)(h')^2+(Z'\circ h)h''.
}
On the other hand $Y$ satisfies \eq{DE+}, $h'=W_{f,g}^{1,0}$ and $h''=(W_{f,g}^{1,0})'=W_{f,g}^{2,0}$ hold on $I$, hence it follows that
\Eq{*}{
(Z''\circ h)\Cdot(W^{1,0}_{f,g})^2+(Z'\circ h)\Cdot W^{2,0}_{f,g}
=\frac{W^{2,0}_{f,g}}{W^{1,0}_{f,g}}\Cdot(Z'\circ h)\Cdot W^{1,0}_{f,g}-\delta(W^{1,0}_{f,g})^2\Cdot(Z\circ h).
}
This reduces to the equality $Z''\circ h=-\delta (Z\circ h)$ on $I$, which, on the interval $h(I)$, is equivalent to
\Eq{*}{
Z''=-\delta Z.
}
Thus, we have proved that $Z:=f\circ h^{-1}$ and $Z:=g\circ h^{-1}$ are solutions to this second-order homogeneous linear differential equation.
The functions $S_{-\delta}$ and $C_{-\delta}$ form a fundamental system of solutions for this differential equation, therefore,
\Eq{*}{
(f\circ h^{-1},g\circ h^{-1})\sim(S_{-\delta},C_{-\delta}).
}
This shows that the relation \eq{gf} is satisfied with $p:=-\delta$, hence, 
from \cor{M3}, we conclude that the assertion (i) holds.

To complete the proof of the theorem it suffices to show that (iv) and (v) are equivalent in the class $\B_2(I)$. If (iv) holds for some $\delta\in\R$, then the differentiability of $W^{1,0}_{f,g}$ implies that $W^{2,1}_{f,g}$ and hence $\Psi_{f,g}$ are differentiable, furthermore,
\Eq{*}{
\frac{\Psi_{f,g}}{(W^{1,0}_{f,g})^2}=-\delta.
}
Differentiating this equation side by side, we obtain
\Eq{*}{
\frac{(W^{1,0}_{f,g})^2\Psi'_{f,g}-2\Psi_{f,g}W^{1,0}_{f,g}W^{2,0}_{f,g}}{(W^{1,0}_{f,g})^4}=0.
}
Simplifying this equality, we can see that (v) must be valid. 

Conversely, if $\Psi_{f,g}$ is differentiable and \eq{Psi} holds, that is, $Y=\Psi_{f,g}$ solves the first-order homogeneous linear differential equation $Y'=2\Phi_{f,g}Y$, then there exists a constant $\delta$ such that 
\Eq{*}{
\Psi_{f,g}
=\delta\exp\bigg(2\int\Phi_{f,g}\bigg)
=\delta\exp\bigg(2\int\frac{(W^{1,0}_{f,g})'}{W^{1,0}_{f,g}}\bigg)
=\delta(W^{1,0}_{f,g})^2,
}
which implies assertion (iv) immediately.
\end{proof}

\Thm{M4}{
Let $(f,g)\in\B_0(I)$. Then the following assertions are equivalent.
\begin{enumerate}[(i)]
\item There exists a continuous strictly monotone function $h:I\to\R$ such that, for all $n\in\N$, $x\in I^n$ and $\lambda\in\Lambda_n$,
\Eq{xl}{
  B_{g,f}(x,\lambda)=A_h(x,\lambda).
}
\item There exists a continuous strictly monotone function $h:I\to\R$ such that, for all $n\in\N$ and $x\in I^n$,
\Eq{x}{
  B_{g,f}(x)=A_h(x).
}
\item There exists a continuous strictly monotone function $h:I\to\R$ and $n\geq3$ such that, for all $x\in I^n$, equation \eq{x} holds.
\item There exists a continuous strictly monotone function $h:I\to\R$ such that equation \eq{xl} holds for all $x\in I^2$ and $\lambda\in\Lambda_2$.
\item There exist $t\in\,]0,\frac12[\,\cup\,]\frac12,1[\,$ and a continuous strictly monotone function $h:I\to\R$ such that equation \eq{xl} holds for all $x\in I^2$ with $\lambda:=(t,1-t)$. 
\item There exist constants $a,b\in\R$ such that 
\Eq{ab+}{
 af+bg=1.
 }
\item Provided that $(f,g)\in\B_2(I)$, $\Psi_{f,g}=0$.
\end{enumerate}
}

\begin{proof}
The implications (i)$\Rightarrow$(ii), (ii)$\Rightarrow$(iii), (i)$\Rightarrow$(iv), and (iv)$\Rightarrow$(v) are obvious. To see that (iii)$\Rightarrow$(v), assume that there exists a continuous strictly monotone function $h:I\to\R$ and $n\geq3$ such that equation \eq{x} is satisfied for all $x\in I^n$. Let $y_1,y_2\in I$ be arbitrary and let $x_1:=y_1$, $(x_2,\dots,x_n):=(y_2,\dots,y_2)$. Applying inequality \eq{x} to the $n$-tuple $x=(x_1,\dots,x_n)$, we get that 
\Eq{*}{
B_{g,f}\big((y_1,y_2),\big(\tfrac1n,\tfrac{n-1}n\big)\big)
=A_h\big((y_1,y_2),\big(\tfrac1n,\tfrac{n-1}n\big)\big)
}
is valid for all $y_1,y_2\in I$. Therefore, assertion (v) holds with $t:=\frac1n$.

To prove the implication (v)$\Rightarrow$(vi), assume that assertion (v) is valid for some continuous strictly monotone function $h$ and $t\in\,]0,\frac12[\,\cup\,]\frac12,1[\,$. Then we have that \eq{B} holds, hence, using \cor{M3}, we get the existence of constants $a,b,c,d\in\R$ with $ad\neq bc$ such that \eq{gf} holds with $p=0$, therefore,  
\Eq{abcd++}{
 1=C_0\circ h=af+bg \qquad\mbox{and}\qquad S_0\circ h=cf+dg.
}
This proves that assertion (vi) is valid. 

Now assume that (vi) holds, i.e., there exist constants $a,b\in\R$ satisfying \eq{ab+}. This equation yields that $a^2+b^2>0$. Define $h:=-bf+ag$. Then we have that $(1,h)\sim(f,g)$, which implies that the Bajraktarevi\'c mean $B_{g,f}$ is identical with the Bajraktarevi\'c mean $B_{h,1}$, which is equal to the quasi-arithmetic mean $A_h$. Therefore, (i) holds, and hence all the assertions from (i) to (vi) are equivalent.
 
To obtain the implication (vi)$\Rightarrow$(vii), assume that $(f,g)\in\B_2(I)$ and that (vi) holds for some constants $a,b\in\R$. Then $af'+bg'=0$ such that $(a,b)\neq(0,0)$. Therefore, $f'$ and $g'$ are linearly dependent. Consequently, we get
\Eq{*}{
W^{2,1}_{f,g}=W^{1,0}_{f',g'}=0.
}
Thus, assertion (vii) is valid.

Finally, it remains to prove the implication (vii)$\Rightarrow$(vi). Let (vii) be satisfied. Then $f$ and $g$ form a system of fundamental solutions of the second-order homogeneous linear differential equation \eq{DE}. In light of assertion (vii), this differential equation reduces to the form
\Eq{*}{
Y''=\Phi_{f,g}Y'.
}
On the other hand, it is clear that $Y=1$ is a solution of this differential equation, therefore it has to be a linear combination of $f$ and $g$. Hence there exist constants $a,b\in\R$ such that \eq{ab+} is satisfied. 
\end{proof}

\def\MR#1{}


\providecommand{\href}[2]{#2}

\end{document}